\documentclass[12pt]{amsart}
\usepackage{enumerate}
\usepackage{float}
\allowdisplaybreaks
\usepackage{tabu}
\usepackage{booktabs}
\usepackage{longtable}
\usepackage{amsfonts,amssymb,amsmath,textcomp}
\usepackage{enumitem}
\usepackage{supertabular}
 2

\usepackage[autostyle]{csquotes}
\theoremstyle{plain}
\newtheorem{thm}{Theorem}[section]

\newtheorem{lem}[thm]{Lemma}
\newtheorem{prop}[thm]{Proposition}

\newcommand{\Gal}[1]{\text{Gal}(#1)}
\def\ord{{\mathrm{ord}}}

\theoremstyle{definition}

\newcommand{\Z}{\mathbb{Z}}
\newcommand{\Q}{\mathbb{Q}}

\begin{document}
\baselineskip 6.1mm

\title
{Primes of higher degree}
\author{Prem Prakash Pandey, Mahesh Kumar Ram}


\address[Prem Prakash Pandey]{Indian Institute of Science Education and Research, Berhampur, India.}
\email{premp@iiserbpr.ac.in}

\address[Mahesh Kumar Ram]{Indian Institute of Science Education and Research, Berhampur, India.}
\email{maheshk@iiserbpr.ac.in}


\date{\today}


\begin{abstract}
Let $K/\Q$ be a cyclic extension of number fields with Galois group $G$. We study the ideal classes of primes $\mathfrak{p}$ of $K$ of residue degree bigger than one in the class group of $K$. In particular, we explore such extensions $K/\Q$ for which there exist an integer $f>1$ such that the ideal classes of primes $\mathfrak{p}$ of $K$ of residue degree $f$ generate the full class group of $K$. It is shown that there are many such fields. These results are used to obtain information on class group of $K$; like rank of $\ell-$torsion of the class group, factors of class number, fields with class group of certain exponents, and even structure of class group in some cases. Moreover, such $f$ can be used to construct annihilators of the class groups.

\end{abstract}

\maketitle{}

\section{Introduction}
For any Galois extension $K/F$ of number fields, we use the notation $G_{K/F}$ for the Galois group of the extension. When $F=\Q$, this will be abbreviated to $G_K$. For any prime ideal $\mathfrak{P}$ of $K$, the degree $[\mathbb{O}_K/\mathfrak{P}: \mathbb{O}_F/ \mathfrak{p}]$ is called the residue degree of the ideal $\mathfrak{P}$ for the extension $K/F$. Here $\mathbb{O}_K, \mathbb{O}_F$ denote the ring of integers of $K$ and $F$ respectively, and $\mathfrak{p}=\mathfrak{P} \cap \mathbb{O}_F$. We shall use the notation $f(\mathfrak{P}|\mathfrak{p})$ for this residue degree. All the prime ideals appearing in the article are assume to be unramified unless stated otherwise explicitly. For any odd prime $\ell$, Kummer showed that the class group of $\Q(\zeta_{\ell})$ is generated by classes of prime ideals of $\Q(\zeta_{\ell})$ of residue degree one. For a general number field $K$, as a consequence of class field theory, one knows that each ideal class in the class groups $C\ell(K)$ of $K$ contains infinitely many prime ideals of residue degree one. In \cite{PPP19}, the first author explored the following problem:\\

{\bf Problem :} Determine extensions $K/F$ for which there exists an integer $f>1$ such that the class group $C\ell(K)$ is generated by prime ideals of residue degree $f$?\\

In relation to this question the following set was introduced. 
\begin{equation*}
\mathbf{R}_{K/F}=\{f \in \mathbb{N}: C\ell(K) \mbox{ is generated by primes of residue degree }f\}.
\end{equation*}
Here the residue degree means the residue degree for the extension $K/F$. The cases where $\mathbf{R}_{K/F}=\{1\}$ are referred as trivial cases.  In \cite{PPP19}, a criteria for some integer $f$ to be in $\mathbf{R}_{K/F}$ was obtained. Under the assumption that $K/F$ is cyclic and $C\ell(F)$ is trivial, it was shown that any $f\in \mathbf{R}_{K/F}$ can be used to construct an annihilator in $\Z[G_{K/F}]$ for the class group of $K$. However, not a single example of extension $K/F$ was shown where $\mathbf{R}_{K/F}$ is non-trivial. In article \cite{NMMRPP}, the authors showed that many cyclotomic and real cyclotomic fields whose class numbers are prime are generated by classes of primes of a fixed residue degree bigger than one. These were the first example of extension $K/F$ with non-trivial $\mathbf{R}_{K/F}$. These were used to obtain some consequences about class number of intermediate fields. \\

This article aims to generalize the works of \cite{PPP19, NMMRPP} in three direction: show that there are many extensions $K/F$ with non-trivial $\mathbf{R}_{K/F}$, obtain the best possible result on annihilator of class groups possible from this study, obtain many strong consequences on class groups of cyclic extensions of $\Q$.  We remark that we obtain these consequences by studying primes of higher degree. For the most part of this paper we use $K$ for a number field which is a cyclic extension of $\Q$ of degree $n$ (with only exception occurring in Theorem \ref{AT} and Section 4). We shall use $C\ell(K)$, $h_K$ to denote the class group of $K$ and the class number of $K$ respectively. We use $H(K)$ to denote the Hilbert class field of $K$. For any group $G$ we use $Aut(G)$ to denote the group of automorphisms of $G$. Our first result in the direction of non-trivial $\mathbf{R}_{K/\Q}$ is the following.

\begin{thm}\label{Main}
Let $K$ be a number field which is cyclic extension of $\Q$ and let $G_K$ be the Galois group of $K/\Q$. Assume that $n$ and $h_K$ are relatively prime. Suppose that there exists a non-trivial element $\sigma \in G_K$ such that any homomorphism $\psi : G_K \longrightarrow Aut(C\ell(K))$ maps $\sigma$ to identity automorphism. Then every divisor of the order of $\sigma$ is in $\mathbf{R}_{K/\Q}$.
\end{thm}


Before proceeding further, we remark that Theorem \ref{Main} generalises the Kummer's result mentioned in the first paragraph, for those cyclic extension $K/\Q$ for which the degree $[K:\Q]$ is relatively prime to the class number $h_K$. Theorem \ref{Main} is quite general result to obtain cyclic extensions $K/\mathbb{Q}$ for which $\mathbf{R}_{K/\Q}$ is non-trivial. If there exists a  divisor $f$ of $n$ which is relatively prime to $|Aut(C\ell(K))|$ then there is an element $\sigma \in G_K$ of order $f$ which satisfies the hypothesis of Theorem \ref{Main}.  We can obtain extensions $K/\mathbb{Q}$ with non-trivial $\mathbf{R}_{K/\Q}$ in more general scenarios as well. In this regard we have the following result.



\begin{thm}\label{main}
Let $K/\Q$ be cyclic extension of number fields and let $G_K$ denote its Galois group. Suppose $C\ell(K) \cong \Z/n_1\Z \oplus \dots \oplus \Z/n_t\Z$, and the class number $h_K$ is relatively prime to the degree $n$ of the extension $K/\Q$. Let $p$ be a prime dividing $n$ with $p^{e_1}$ being the highest power of $p$ dividing $n$, and let $p^{e_2}$ be the highest power of $p$ in $|Aut(C\ell(K))|$ with $e_2 \geq 0$. If $e_1>e_2$ then $p^i \in \mathbf{R}_{K/\Q}$ for all $i \in \{0, 1, \dots, e_1-e_2\}$.
\end{thm}

Before proceeding further, we remark that above results help in exhibiting extensions $K/\Q$ for which $\mathbf{R}_{K/\Q}$ is non-trivial. For such extensions, most of the times, it is possible to choose intermediate fields $F$ such that $\mathbf{R}_{K/F}$ is non-trivial. \\

Now we turn to the second goal of the article, that is, construction of annihilators of class groups. Annihilators of class groups are important objects in algebraic number theory \cite{FT88, JT81, YFB04, RS08,  BBM14, JSBT18,CGRC21}. It was established in \cite{PPP19} that every element in $\mathbf{R}_{K/F}$ can be used to construct annihilators of $C\ell(K)$ provided the Galois group $G_{K/F}$ is cyclic and the class number $h_F=1$. The annihilators of \cite{PPP19} were extended to more general extensions $K/F$ in \cite{NMMRPP}.The method of \cite{NMMRPP} yields annihilator corresponding to non-trivial element in $\mathbf{R}_{K/F}$ even for non-abelian extensions. In this regard we obtain a stronger result than mentioned in \cite{NMMRPP}. For this part of discussion, we emphasize that $K/F$ is any extension of number fields (may be non-abelian). Let $f \in \mathbf{R}_{K/F}$, and $H_f^1, \dots, H_f^r$ denote all the cyclic subgroups of $G_{K/F}$ of order $f$. Put $H_f=\cap_{i=1}^r H_f^i$ and let $T=\{\sigma_1, \dots, \sigma_t\}$ denote a set of coset representatives of $G_{K/F}/H_f$. Define
\begin{equation}\label{ann}
\theta_f=\sum_{i=1}^t \sigma_i.
\end{equation}
The following theorem generalizes Theorem 1.2 of \cite{NMMRPP}.
\begin{thm}\label{AT}
If the class group $C\ell(F)$ is trivial then the element $\theta_f \in \Z[G_{K/F}]$ annihilates the class group $C\ell(K)$.
\end{thm}

Lastly we mention some consequences of study of primes of higher degree, on the class groups. For example, Kummer wanted to prove that the cubic subfield of $\Q(\zeta_{163})$ can not have class number $2$. Using the $\Z[G_{\Q(\zeta_{163})}]$ structure of the class group of $\Q(\zeta_{163})$ this can be proved (see \cite{CG21}). In a private communication Prof. Greither informed us that this can be proved for many other cyclic extensions. We prove a much more general result in this direction. For any integer $m$, let $\phi(m)$ denote the value of Euler-totient function at $m$. 

\begin{thm}\label{A1}
Let $n$ and $m$ be two relatively prime integers bigger than $1$. Assume that $n$ is relatively prime to $\phi(m)$.There is no cyclic extension $K/\Q$ of degree $n$ whose class group is cyclic of order $m$. 
\end{thm}

Theorem \ref{A1} shows that there is some constraint on prime factors of the class number of a cyclic extension $K/\Q$. On the other hand, if we know factors of class number, primes of higher degree can be used to obtain the rank of class groups. For this let $n=q_1^{e_1} \ldots q_r^{e_r}$ be the factorization of $n$ as product of primes. For any prime number $\ell$ not dividing $n$, let $s_i$ denote the multiplicative order of $\ell$ modulo $q_i$. Put $s=\min \{s_1, \ldots, s_r\}$. Let $C\ell(K)[\ell]$ denote the $\ell$-torsion of the class group $C\ell(K)$.  That is,
$$C\ell(K)[\ell]=\{[\mathfrak{a}] \in C\ell(K): [\mathfrak{a}]^{\ell}=1\}.$$
Then we prove the following theorem.
\begin{thm}\label{A11}
Let $K/\Q$ be cyclic extension of degree $n=q_1^{e_1} \ldots q_r^{e_r}$ and $\ell$ be any prime number not dividing $n$. If $\ell$ divides $h_K$ then $C\ell(K)[\ell]$ has a subgroup isomorphic to $(\Z/\ell\Z)^s$.
\end{thm}

We remark that the $\ell-$torsion subgroup $C\ell(K)[\ell]$ are very important (see \cite{HV06, EV07, PCW20, PCW21}) and a lot of results are obtained on them in last twenty years or so (see \cite{LB05, HV06, EV07, EPW17, BSTTTZ20, PCW20, WJ21, KP22}). Our result is purely based on the study of primes of higher degree. Lastly we mention one more result, on class groups, obtained from our study of primes of higher degree. 
\begin{thm}\label{A3}
Consider the cyclotomic extension $\Q(\zeta_{\ell})$ for some prime $\ell$. Assume that $(\ell-1)$ is relatively prime to $C\ell(\Q(\zeta_{\ell}))$. Let $f$ be a divisor of $(\ell-1)$ which is relatively prime to both $(\ell-1)/f$ and $|Aut(C\ell(\Q(\zeta_{\ell})))|$. Then the exponent of the class group of the subfield of degree $f$ is at most $(\ell-1)/f$.
\end{thm}

In Section 2 we recall some results from group theory and some results from class field theory which will be useful. Section 3 contains a proof of theorem \ref{Main} and Theorem \ref{main}.  Section 4 contains some more extensions $K/\Q$ (non-cyclic) with non-trivial $\mathbf{R}_{K/\Q}$. In Section 5 we prove Theorem \ref{AT}.  Sections 6 contains a proof of the consequences Theorem \ref{A1}, Theorem \ref{A11} and Theorem \ref{A3}. In Section 6 we mention some more consequences of the study of primes of higher degree. The applications mentioned here are not exhaustive and some more such applications will appear in the second author's thesis. The last section ends with providing explicit examples of extensions $K/\Q$ for which $\mathbf{R}_{K/\Q}$ is non-trivial. The second author's thesis will list many more examples produced from this method.

\vspace{5mm}
\section{Preliminaries}
We begin with a Galois extension $L/F$ of number fields. For any unramified prime ideal $\mathfrak{P}$ of $L$, the Frobenius automorphism of $\mathfrak{P}$ is the unique element $\sigma$ of the decomposition group $D_{\mathfrak{P}}$ satisfying 
$$\sigma (x) \equiv x^{N(\mathfrak{p})} \pmod {\mathfrak{P}}.$$
Here $\mathfrak{p}=\mathfrak{P} \cap \mathbf{O}_F$ and $N(\mathfrak{p})$ is the norm for the extension $F/\Q$. The Frobenius automorphism of $\mathfrak{P}$ will be denoted by $\left(\frac{\mathfrak{P}}{L/F}\right)$, and by $\left(\frac{\mathfrak{p}}{L/F}\right)$ in case $G_{L/F}$ is abelian. We recall the following result from \cite{GJ96}.
\begin{lem}\label{FR}
Let $L/F$ be a Galois extension of number fields and $E$ be an intermediate field. For any unramified prime ideal $\mathfrak{P}$ of $L$, we have
$$\left( \frac{\mathfrak{P}}{L/F} \right)^{f(\mathbf{p}/\mathfrak{p})}=\left( \frac{\mathfrak{P}}{L/E}\right).$$
Here $\mathbf{p}$ and $\mathfrak{p}$ are the prime ideals below $\mathfrak{P}$ in $E$ and $F$ respectively.
\end{lem}
Let $H(L)$ denote the Hilbert class field of the number field $L$. The following result is an important consequences of class field theory \cite{LW91}
\begin{lem}\label{HCF}
The Galois group $G_{H(L)/L}$ is isomorphic to the class group $C\ell(L)$. The isomorphism is induced by the Artin map
$$[\mathfrak{P}] \longmapsto \left( \frac{\mathfrak{P}}{L/F} \right).$$
\end{lem}
The groups $G_{H(L)/L}$ and $C\ell(L)$ can be used in place of each other, up to isomorphism, and the justification is offered by Lemma \ref{HCF}. Next, we recall the $\check{C}$ebotarev density Theorem. For any $\sigma \in G_{L/F}$ we consider the set $P_{L/F}(\sigma)$ of prime ideals $\mathfrak{p}$ of $F$ such that there is a prime $\mathfrak{P}$ of $L$ above $\mathfrak{p}$ satisfying
$$\sigma =\left( \frac{\mathfrak{P}}{L/F} \right),$$
also we use $C_{\sigma}$ to denote the conjugacy class of $\sigma$ in $G_{L/F}$. Now we state the $\check{C}$ebotarev density theorem \cite{JN99}.
\begin{thm}\label{CDT}
Let $L/F$ be a Galois extension with Galois group $G_{L/F}$. For every $\sigma \in G_{L/F}$, the density of the set $P_{L/F}(\sigma)$ is positive and equals to $\frac{|C_{\sigma}|}{[L:F]}$. 
\end{thm}
The next lemma can be proved elementarily \cite{PPP19}.
\begin{lem}
The extension $H(L)/F$ is Galois.
\end{lem}
For any Galois extension of number fields $L/F$, there is following exact sequence
$$1 \longrightarrow G_{H(L)/L} \longrightarrow G_{H(L)/F} \longrightarrow G_{L/F} \longrightarrow 1.$$
We need the following theorem of Wyman \cite{BW73} (also see \cite{GCMR88}).
\begin{thm}\label{splitting}
If $K/\Q$ is a cyclic extension then the sequence
$$1 \longrightarrow G_{H(K)/K} \longrightarrow G_{H(K)} \longrightarrow G_K \longrightarrow 1$$
splits. Consequently, there exists a homomorphism $\psi : G_K \longrightarrow Aut(G_{H(K)/K})$ such that
$$G_{H(K)} \cong G_{H(K)/K} \rtimes_{\psi} G_K.$$
\end{thm}

Now we state some results on groups. The following lemma can be easily proved.
\begin{lem}\label{el}
Let $H$ be an abelian group of order $n$. Assume that $n=r.s$ for relatively prime integers $r$ and $s$. If $H_1$ and $H_2$ are subgroups of orders $r$ and $s$ respectively then the set $ H_2$ is a complete set of coset representatives of $H/H_1$.
\end{lem}
Now we recall some results on the automorphism group of finite abelian groups (see \cite{CHDR07, JP04, AR1907}). We begin with the following result on automorphism group of finite abelian $p-$groups. Let $H$ be a finite abelian $p-$group and
\begin{equation}\label{41}
H \cong \Z/p^{e_1}\Z \oplus \dots \oplus \Z/p^{e_t}\Z.
\end{equation}
Define $$d_i=\max \{j: e_j=e_i\} \mbox{ and } c_i=\min \{j: e_j=e_i\}.$$
\begin{prop}\label{CAAG}
Let $H$ be as in (\ref{41}) and $Aut(H)$ denote the group of automorphisms of $H$. Then
\begin{equation}\label{42}
|Aut(H)|=\prod_{i=1}^t \left(p^{d_i}-p^{(i-1)}\right) \prod_{i=1}^t \left(p^{e_i}\right)^{(t-d_i)} \prod_{i=1}^t \left( p^{(e_i-1)}\right)^{(t-c_i+1)}.
\end{equation}
\end{prop}
The next lemma is elementary and facilitates us with the automorphism group of any abelian group.
\begin{lem}\label{ARPOP}
Let $G_1$ and $G_2$ be finite abelian groups of relatively prime orders. Then 
$$Aut(G_1 \oplus G_2) \cong Aut(G_1) \oplus Aut(G_2).$$
\end{lem}
Next, we prove an elementary lemma about the order of certain elements in a semidirect product of groups.
\begin{lem}\label{SPL}
Let $G_1, G_2$ be two finite groups and $\psi :G_2 \longrightarrow Aut(G_1)$ be a group homomorphism. Let $g_1 \in G_1, g_2 \in G_2$ be of orders $n_1$ and $n_2$ respectively. If $\psi(g_2)$ is identity automorphism then the element $(g_1, g_2) $ in the semidirect product $G_1 \rtimes_{\psi} G_2$ has order $lcm~ (n_1, n_2)$.
\end{lem}
\begin{proof}
Since $\psi(g_2)$ is identity automorphism of $G_1$, we have $\psi(g_2)(g)=g$ for all $g\in G_1$. In particular, $\psi(g_2)(g_1^i)=g_1^i$ for any integer $i$. As a result, we have
$$(g_1,g_2)^2=(g_1 \psi(g_2)(g_1), g_2^2)=(g_1^2,g_2^2).$$
Proceeding inductively yields $(g_1, g_2)^n=(g_1^n, g_2^n)$.  The Lemma follows at once.
\end{proof}

\vspace{5mm}
\section{Proofs of Theorem \ref{Main} and Theorem \ref{main}}
\begin{proof} (Theorem \ref{Main}) As $K/\Q$ is a cyclic extension, from Theorem \ref{splitting} it follows that there is a homomorphism $\psi : G_K \longrightarrow Aut(G_{H(K)/K})$ inducing an isomorphism
\begin{equation}\label{MTE1}
\Psi : G_{H(K)} \longrightarrow G_{H(K)/K} \rtimes_{\psi} G_K.
\end{equation}
Since $C\ell(K)$ is an abelian group and $C\ell(K) \cong G_{H(K)/K}$, we have
$$G_{H(K)/K} \cong \Z/n_1\Z \oplus \dots \oplus \Z/n_t\Z$$ for some integers $n_1, \dots, n_t$. Let $\sigma_1, \dots, \sigma_t \in G_{H(K)/K}$ be such that order of $\sigma_i$ is $n_i$ for each $i$ and 
\begin{equation*}
G_{H(K)/K}=\left\{\prod_{i=1}^t \sigma_i^{j_i}: 1 \leq j_i \leq n_i , \forall i\right\}.
\end{equation*}
We are given that $\psi(\sigma)=Id$. Now suppose the order of $\sigma$ equals to $f$.  Then $f$ is relatively prime to each $n_i$, hence we have
\begin{equation}\label{MTE2}
G_{H(K)/K}=\left\{\prod_{i=1}^t \sigma_i^{fj_i}: 1 \leq j_i \leq n_i , \forall i\right\}.
\end{equation}
From Lemma \ref{SPL}, the element $(\sigma_i, \sigma)$ have orders $n_i.f$ for each $i$.  Because of the isomorphism in (\ref{MTE1}), there are elements $\tau_i \in G_{H(K)}$ such that $\Psi(\tau_i)=(\sigma_i, \sigma)$. Applying the $\check{C}$ebotarev density theorem we obtain unramified prime ideals $\mathfrak{P}_i$ of the field $H(K)$ such that
\begin{equation}\label{MTE3}
\left(\frac{\mathfrak{P}_i}{H(K)/\Q} \right)= \tau_i, \mbox{ for each }i=1, \dots, t.
\end{equation} 
We let $\mathfrak{p}_i$ and $p_i$ denote the primes below $\mathfrak{P}_i$ of the fields $K$ and $\Q$ respectively. We have the following identities
$$f(\mathfrak{P}_i|p_i)= \ord{~\tau_i}= n_i.f \mbox{ and } f(\mathfrak{P}_i|p_i)=f(\mathfrak{P}_i|\mathfrak{p}_i).f(\mathfrak{p}_i|p_i).$$
Note that $f(\mathfrak{P}_i|\mathfrak{p}_i)$ is a divisor of $[H(K):K]$ and $f(\mathfrak{p}_i|p_i)$ is a divisor of $n$. From this we obtain the following
\begin{equation}\label{MTE4}
f(\mathfrak{P}_i|\mathfrak{p}_i)= n_i \mbox{ and } f(\mathfrak{p}_i|p_i)=f.
\end{equation}
Put $$\left(\frac{\mathfrak{P}_i}{H(K)/K} \right)= \tau_i^{'}, \mbox{ for each }i=1, \dots, t.$$
Then from Lemma \ref{FR}, we get $\tau_i^f=\tau_i^{'}$ for each $i$. Thus $\ord{~\tau_i'}=n_i$. We claim that the Galois group $G_{H(K)/K}$ is generated by $\tau_1', \dots, \tau_t'$.\\
\noindent We have $$\Psi(\tau_i'^{j_i})=\Psi(\tau_i^{fj_i})=(\sigma_i, \sigma)^{fj_i}.$$ Using Lemma \ref{SPL}, we see that
\begin{equation}\label{MTE5}
\Psi(\tau_i'^{j_i})=(\sigma_i^{fj_i}, 1).
\end{equation} 
Since $G_{H(K)/K}$ embeds inside $G_{H(K)/K} \rtimes_{\psi} Gal(K/Q)$ via the map $g \mapsto (g,1)$ and $\Psi$ is an isomorphism, our claim follows from (\ref{MTE2}) and (\ref{MTE5}).\\
As the extension $H(K)/K$ is abelian we have
$$\left(\frac{\mathfrak{p}_i}{H(K)/K} \right)=\left(\frac{\mathfrak{P}_i}{H(K)/K} \right)=\tau_i'.$$
Using Lemma \ref{HCF}, we conclude that $C\ell(K)$ is generated by the classes of prime ideals $\mathfrak{p}_i$. This shows that $f \in \mathbf{R}_{K/\Q}$.\\
For any divisor $f'$ of $f$, working with $\sigma^{f/f'}$ in place of $\sigma$, we obtain $f' \in \mathbf{R}_{K/\Q}$. This finishes the proof of Theorem \ref{Main}.\\
\end{proof}

To generate examples of number fields $K$ with non-trivial $\mathbf{R}_{K/\Q}$ we state a weaker version of Theorem \ref{Main} which can be readily applied. This can be easily derived from Theorem \ref{Main}.
\begin{thm}\label{CMain}
Let $K/\Q$ be cyclic with Galois group $G_K$. Assume that $n$ and $h_K$ are relatively prime. Let $f$ be a divisor of $n$ which is relatively prime to $|Aut(C\ell(K))|$. Then every divisor of $f$ is in $\mathbf{R}_{K/\Q}$.
\end{thm}

Along the lines of the proof of Theorem \ref{Main}, we can prove the following generalization of Theorem \ref{Main}.
\begin{thm}\label{GMain}
Let $K/\Q$ be cyclic with Galois group $G_K$. Let $$C\ell(K) \cong \Z/n_1\Z \oplus \dots \oplus \Z/n_t\Z$$ be the decomposition of the class group into the invariant factors with $n_1|n_2| \dots |n_t$. Let $\sigma \in G_K$ be a non-trivial element of order $f$ which is relatively prime to $n_t$ and\\
(1) for any homomorphism $\psi : G_K \longrightarrow Aut(H(K))$ we have $\psi(\sigma)=Id$,\\
(2) the factor $\frac{n}{f}=q$ is a prime relatively prime to $f$.\\
Then the subgroup of $C\ell(K)$ generated by prime ideals of residue degree $f$ contains a subgroup isomorphic to $$\Z/n_1^{'}\Z \oplus \dots \oplus \Z/n_t^{'}\Z,$$ where $n_i^{'}=\frac{n_i}{q^{r_i}}$, with $q^{r_i}$ being the highest power of $q$ dividing $n_i$.
\end{thm}

One can obtain more generalizations, some such generalizations will appear in second author's thesis. \\

Now we give the proof of Theorem \ref{main}.

\begin{proof}(Theorem \ref{main})
The Galois group $G_K$ is cyclic and $p^{e_1}$ divides $n$, so there exist an element $\sigma \in G_K$ of order $p^{e_1}$. For any homomorphism $\psi : G_K \longrightarrow Aut(C\ell(K))$, the order of the element $\psi(\sigma)$ is a divisor of $p^{e_1}$. Since $p^{e_2}$ is the highest power of $p$ which divides $|Aut(C\ell(K))|$, we conclude that the order of $\psi(\sigma)$ is at most $p^{e_2}$. In particular,
$$\psi(\sigma)^{p^{e_2}}=Id.$$
Thus for $\sigma'=\sigma^{p^{e_2}}$, we have
$$\psi(\sigma')=\psi(\sigma^{p^{e_2}})=\psi(\sigma)^{p^{e_2}}=Id.$$ 
Since the order of $\sigma'$ is $p^{e_1-e_2}$, the theorem follows from Theorem \ref{Main}.
\end{proof}
\vspace{5mm}

\section{More families of extensions $K/\Q$ with non-trivial $\mathbf{R}_{K/\Q}$}

In this section we do not require the extension $K/\Q$ to be cyclic. We need the following analogue of Theorem \ref{splitting}.
\begin{prop}\label{Splitting-4}
Let $K/\Q$ be a Galois extension of number fields of prime power degree, say $n=q^t$ for some prime $q$. Assume that $q$ does not divide the class number $h_K$. Then the sequence
\begin{equation}\label{43}
1 \longrightarrow C\ell(K) \longrightarrow G_{H(K)} \longrightarrow G_K \longrightarrow 1
\end{equation}
splits. In particular, we have 
\begin{equation}\label{44}
G_{H(K)} \cong C\ell(K) \rtimes G_K.
\end{equation}
\end{prop}
\begin{proof}
It is enough to find a subfield $L$ of $H(K)$ such that $K \cap L=\Q$ and $H(K)=KL$. Let $T$ be a Sylow $q-$subgroup of $G_{H(K)}$ and $L$ be the fixed field of $T$. As $[H(K):K]=h_K$ is relatively prime to $q$, we see that
$$[H(K):L]=|T|=q^t.$$
From this, we immediately obtain $$[L:\Q]=\frac{[H(K):\Q]}{q^t}=\frac{[H(K):\Q]}{n}=h_K.$$
Since $[L:\Q]$ and $n$ are relatively prime, we obtain at once
$$K \cap L=\Q \mbox{ and } [KL:\Q]=n[L:\Q]=[H(K):\Q].$$
Now the proposition follows.

\end{proof}
Now we are in a position to prove the main theorem of this section.
\begin{thm}\label{T41}
Let $K/\Q$ be a Galois extension of degree $q^t$ for some $t>1$ and an odd prime $q$. Assume that the class number factors as $h_K=p_1^{e_1} \dots p_u^{e^u}$ with distinct primes $p_i$. If $q$ does not divide $|Aut(C\ell(K))|$ and is relatively prime to $h_K$ then $q \in \mathbf{R}_{K/\Q}$.
\end{thm}
\begin{proof}
From Proposition \ref{Splitting-4} it follows that there is a homomorphism $$\psi: G_K \longrightarrow Aut(C\ell(K))$$ inducing an isomorphism $$\Psi: G_{H(K)} \longrightarrow C\ell(K) \rtimes_{\psi} G_K.$$
Let $\sigma \in G_K$ be an element of order $q$.  As $q$ does not divide $|Aut(C\ell(K))|$, we must have $\psi(\sigma)=Id$.\\
Now rest of the proof proceeds along the same line as that of Theorem \ref{Main}.
\end{proof}

Before stating our next Theorem, we recall the following lemma from elementary group theory.
 
\begin{lem}\label{L43}
Let $n=q_1^{e_1}q_2^{e_2}\cdots q_r^{e_r}$ be an integer, where $q_i$'s are distinct odd primes, and $e_i\text{'s}\geq 1$ for all $i\in\{1,2,\dots,r\}$. Then
\begin{enumerate}[label=(\roman*)]
    \item $Aut(\frac{\mathbb{Z}}{n\mathbb{Z}})\cong \prod\limits_{i=1}^r\frac{\mathbb{Z}}{q_i^{e_i-1}(q_i-1)\mathbb{Z}}$
    
    \item the number of elements of order $2$ in $Aut(\frac{\Z}{n\Z})$ is $2^r-1$.
\end{enumerate}
\end{lem}
Another lemma along the same lines is the following.
\begin{lem}\label{L44}
Let $G_1$ and $G_2$ be finite abelian groups such that $G_1\cong (\frac{\Z}{2\Z})^{t_1}$, and the total number of elements of order $2$ in the group $G_2$ is $2^{t_2}$, where $t_1$ and $t_2$ are positive integers such that $t_1>t_2\geq 1$. Then, for any homomorphism $\psi: G_1\rightarrow G_2$, the kernel of $\psi$ has an element of order $2$.
\end{lem}
\begin{proof}
For any element $ \sigma$ of  $(\frac{\Z}{2\Z})^{t_1}$, order of $\psi(\sigma)$ is either $1$ or $2$.  Since $G_2$ has $2^{t_2}$ elements of order $2$,  and $2^{t_1}>2^{t_2}+1$, it follows that there is a non-trivial element in the kernel of $\psi$. 
\end{proof}
As an analogue to Theorem \ref{T41} for the case $q=2$ we prove the following theorem.
\begin{thm}\label{T45}
Let $K/\mathbb{Q}$ be a Galois extension of degree $2^{t_1}$ such that $G_K \cong (\frac{\Z}{2\Z})^{t_1}$, where $t_1\geq 2$ is an integer. Assume that the class group $C\ell(K)$ of $K$ is cyclic of odd order, and let $h_k=p_1^{e_1} \dots  p_{t_2}^{e_{t_2}}$. If $t_1 >t_2$, then $2\in\mathbf{R}_{K/\Q}.$
\end{thm}
\begin{proof}
From Proposition \ref{Splitting-4} it follows that there is a homomorphism $$\psi: G_K \longrightarrow Aut(C\ell(K))$$ inducing an isomorphism $$\Psi: G_{H(K)} \longrightarrow C\ell(K) \rtimes_{\psi} G_K.$$
From Lemma \ref{L43} and Lemma \ref{L44}, we find that there exists an element $\sigma$ in $G_K$ of order $2$ such that $\psi(\sigma)=Id$.\\
Now rest of the proof proceeds along the same line as that of Theorem \ref{Main}.
\end{proof}

\vspace{5mm}
\section{Proof of Theorem \ref{AT}}
We need the following result about groups.
\begin{lem}\label{L2}
Let $H_2$ be a subgroup of $H_1$ and $H_1$ be a subgroup of a finite group $H$. If $T$ is a set of representatives of all the left cosets of $H_2$ in $H$, then there exist sets $T_1$ and $T_2$ of representatives of all the left cosets of $H_2$ in $H_1$ and $H_1$ in $H$ respectively such that for each $g \in T$ there is a unique pair $(g', g^{''}) \in T_1 \times T_2$ and an element $h \in H_2$ satisfying
\begin{equation}\label{L2e1}
g=g^{''}g'h
\end{equation}
and vice-versa.
\end{lem}
\begin{proof}
Let $T=\{z_1, \dots, z_t\}$. Since $T$ is a set of representatives of the set of left cosets of $H_2$ in $H$,
\begin{equation}\label{Pe11}
H= \bigsqcup_{i=1}^t z_iH_2.
\end{equation}
Let $T_1$ be the set of all $z_i's$ such that $z_iH_2 \subset H_1$. We consider the largest subset $T_2' \subset T$ with the property
 $$z_i \not \in H_1 \mbox{ holds for all but one } z_i \in T_2'.$$
Let $T_2$ be the largest subset of $T_2'$ such that for any two distinct elements $z_i, z_j \in T_2$ we have $z_iH_1 \not = z_jH_1$. \\
From the definition of $T_1$ and (\ref{Pe11}) it follows that $T_1$ is a set of representatives for the set of left cosets of $H_2$ in $H_1$. That is,
\begin{equation}\label{Pe12}
H_1=\bigsqcup_{z_i \in T_1}z_iH_2.
\end{equation}
By definition, there is exactly one $y_i\in T_2'$ such that $y_iH_1=H_1$. Clearly the element $y_i \in T_2$. Let $g\in H$ be such that $gH_1 \not = H_1$. Now, $gH_2=z_iH_2$ for some $z_i \in T$. As $H_2 \subset H_1$ and $gH_1 \neq H_1$, we must have  $z_iH_1 \neq H_1$. Thus $z_i \in T_2'$. From the definition of $T_2$, there exists unique $z_j \in T_2$ such that $z_iH_1=z_jH_1$. It is readily seen that $gH_1=z_jH_1$. From this, we conclude that
\begin{equation}\label{Pe13}
H=\bigsqcup_{z_j \in T_2}z_jH_1.
\end{equation}
From (\ref{Pe12}) and (\ref{Pe13}) it follows that
\begin{equation}\label{Pe14}
H=\bigsqcup_{\substack{z_i \in T_1}\ z_j \in T_2}z_jz_iH_2.
\end{equation}
Now, the assertion (\ref{L2e1}) follows from (\ref{Pe11}) and (\ref{Pe14}). This proves the lemma.
\end{proof}

\begin{proof}(Proof of Theorem \ref{AT})
Let $\wp$ be an unramified prime ideal of $L$ of residue degree $f$. We let $D_f$ denote the decomposition subgroup of $G_{L/F}$ for the prime ideal $\wp$.  Let $H_f^1, \ldots, H_f^r$ denote all the cyclic subgroups $G_{L/F}$ of order $f$. Then $D_f=H_f^i$ for some $i$. If we put $H_f=\bigcap_{i=1}^r H_f^i$, then there exist a set $T=\{\sigma_1, \dots, \sigma_t\}$ of representatives of left cosets of $H_f$ in $G_{L/F}$ such that
\begin{equation}\label{P1e1}
\theta_f=\sum_{\sigma_i \in T} \sigma_i.
\end{equation}
From Lemma \ref{L2}, there exist subsets $T_1, T_2$ of $T$ such that $T_1$ is a set of representatives of left cosets of $H_f$ in $D_f$ and $T_2$ is a set of representatives of left cosets of $D_f$ in $G_{L/F}$. Moreover, there are elements $\tau_{ij}$ in $H_f$ such that 
\begin{equation}\label{P1e2}
\theta_f=\sum_{\sigma_i' \in T_1, \sigma_j^{''} \in T_2,} \sigma_j^{''} \sigma_i' \tau_{ij}.
\end{equation}
As a result,
\begin{equation}\label{P1e3}
\wp^{\theta_f}=\prod_{\sigma_i' \in T_1, \sigma_j^{''} \in T_2,} \sigma_j^{''} \sigma_i' \tau_{ij}(\wp).
\end{equation}
Since $H_f \subset D_f$, we have $\tau_{ij}(\wp)=\wp$. Consequently.
\begin{equation}\label{P1e4}
\wp^{\theta_f}=\prod_{\sigma_j^{''} \in T_2} \prod_{\sigma_i' \in T_1}, \sigma_j^{''}  \sigma_i'  (\wp).
\end{equation}
For each $\sigma_i' \in T_1$, we have $\sigma_i(\wp)=\wp$. As a result, we obtain
\begin{equation}\label{P1e5}
\wp^{\theta_f}=\prod_{\sigma_j^{''} \in T_2}  \sigma_j^{''}(\wp^{t_1}) \mbox{ here } t_1=|T_1|.
\end{equation}
If $\mathfrak{p}$ is the prime ideal of $F$ lying below $\wp$ then, from the factorization theorem of Dedekind, we see that
\begin{equation}\label{P1e6}
\mathfrak{p} \mathcal{O}_L=\prod_{\sigma_j^{''} \in T_2} \sigma_j^{''}(\wp).
\end{equation}
From equations (\ref{P1e5}) and (\ref{P1e6}) we obtain
$$\wp^{\theta_f}=(\mathfrak{p}\mathcal{O}_L)^{t_1}.$$
Since $h_F=1$, the ideals $\mathfrak{p}$ and $\mathfrak{p}\mathcal{O}_L$ must be principal. From this we conclude that $\wp^{\theta_f}$ is principal. As the class group $C\ell(L)$ is generated by unramified prime ideals of residue degree $f$, the theorem follows.
\end{proof}

\vspace{5mm}

\section{applications}
In this section we prove Theorem \ref{A1}, Theorem \ref{A11} and Theorem \ref{A3}.  We begin with the proof of Theorem \ref{A1}.
\begin{proof}(Theorem \ref{A1})
If possible let $K/\Q$ be a cyclic extension such that the class group $C\ell(K)$ is cyclic of order $m$. Then $|Aut(C\ell(K))|=\phi(m)$.  From Theorem \ref{CMain}, we see that $n \in \mathbf{R}_{K/\Q}$. This means the class group $C\ell(K)$ is generated by prime ideals of residue degree $n$.  As $n=n$, prime ideals of $K$ of residue degree $n$ are principal. Consequently the class group $C\ell(K)$ must be trivial. This proves Theorem \ref{A1}.
\end{proof}

Now we prove Theorem \ref{A11}. We prove Theorem \ref{A11} for $n=q$ an odd prime.  With this assumption the arguments can be easily conveyed.  Let $s$ be the multiplicative order of $\ell$ modulo $q$.  Recall that $G_{H(K)/K} \cong H_K$. Suppose $h_k=\ell^t m$ with $\ell \nmid m, t>0$. 

Let $H$ denote the $\ell-$ Syllow subgroup of $C\ell(K)$, then there are integers $n_1, \ldots, n_r$ such that
\begin{equation}\label{ea1}
H \cong (\Z/\ell^{n_1}\Z) \oplus \ldots \oplus (\Z/\ell^{n_r}\Z).
\end{equation}
Next, let $M$ be the unique subgroup of $C\ell(K)$ of order $m$ and let $L$ denote the fixed field of $M$ for the extension $H(K)/K$. Then $L/K$ is Galois and $G_L/K \cong H$.
 \begin{prop}\label{Pa3}
The extension $L/\Q$ is Galois.
\end{prop}
\begin{proof}
Suppose the extension $L/\Q$ is not Galois. This means that there exists an embedding $\sigma: L\rightarrow \Bar{\Q}$ such that $\sigma(L)\neq L$. Put $L_\sigma=\sigma(L)$. Clearly, we have  $[L: \Q]=[L_\sigma: \Q]=\ell^t[K:\Q]$ and $[L: K]=[L_\sigma: K]=\ell^t$.  From this it follows that 
$$[L\cap L_\sigma :K]= \ell^i \;\text{for some}\; 0\leq i< t.$$
Since $L/K$ is Galois, it follows that
$$[L L_\sigma :K]= \ell^{(2t-i)}.$$
From this we get that $\ell^{(t+1)}|h_K$. This contradiction proves the proposition.
\end{proof}
Now, we have the following short exact sequence\\
\begin{equation}\label{ea2}
1 \longrightarrow G_{L/K} \longrightarrow G_L \longrightarrow G_K \longrightarrow 1.
\end{equation}
Analogous to Theorem \ref{splitting}, the exact sequence (\ref{ea2}) splits. Consequently, one has
\begin{equation}\label{ea21}
G_L \cong G_{L/K} \rtimes G_K.
\end{equation}
From (\ref{ea1}), it follows that $$C\ell(K)[\ell] \cong \left(\Z/\ell\Z\right)^r. $$
We need to show that $r\geq s$. Suppose this is not true. From Proposition \ref{CAAG} it follows that $|Aut(G_{L/K})|$ is relatively prime to $q$. Thus if $\sigma$ is a generator of $G_K$, then for any homomorphism $\psi : G_K \longrightarrow Aut(G_{L/K})$ we see that $\psi(\sigma)=Id$. Now we consider the following set
\begin{equation*}\label{e3}
\mathbf{R}_{K/\Q}[\ell]=\{f \in \mathbb{N}: C\ell(K)[\ell] \mbox{ is generated by primes of residue degree }f\}.
\end{equation*}
Note that $n \in \mathbf{R}_{K/\Q}[\ell]$ implies that $C\ell(K)[\ell]$ is trivial. Thus, showing $n \in \mathbf{R}_{K/\Q}[\ell]$ will complete the proof. This follows from the next theorem.
\begin{thm}
Let $\tau \in G_K$ be an element of order $f$ and assume that for any homomorphism $\psi: G_K\rightarrow \text{Aut}(G_{L/K})$ we have $\psi(\tau)=Id$. Then each divisor of $f$ is an element of $\mathbf{R}_{K/\Q}[\ell]$.
\end{thm}
\begin{proof}
Let $L$ be the field, defined as above, such that $G_{L/K} \cong H $.  From Proposition \ref{Pa3} it follows that the extension $L/\Q$ is a Galois extension. From (\ref{ea21}) we obtain a homomorphism $\psi :G_K   \rightarrow\text{Aut}(G_{L/K})$ inducing an isomorphism 
$$\Psi: G_L \rightarrow G_{L/K }\rtimes_{\psi} G_K.$$ 
We choose elements $\sigma_1, \sigma_2, \dots,\sigma_r$ form $G_{L/K}$ such that the order of $\sigma_i$ is $\ell^{n_i}$ for all $i\in \{1,2,\dots, r\}$ and 
$$G_{L/K}=\{\prod\limits_{i=1}^{r}\sigma_i^{j_i}: j_i\in \Z \;\text{for each}\;i\}.$$ 
As $f$ and $q$ are relatively prime, it follows that
\begin{equation}\label{ea4}
G_{L/K}= \{\prod\limits_{i=1}^{r}\sigma_i^{fj_i}: j_i\in \Z \;\text{for each}\;i\}.
\end{equation}
Given that $\psi(\tau)=Id$, we see that the order of $(\sigma_i, \tau)$ in $G{L/K}\rtimes_{\psi} G_K$ is $\ell^{n_i}f$.

Since $\Psi$ is an isomorphism, there exist elements $\tau_i$'s in $G_L$ of order $\ell^{n_i}f$ such that $\Psi(\tau_i)= (\sigma_i, \tau)$. Applying $\Check{C}$ebotarev density theorem for each $\tau_i$, we obtain prime ideals $\mathfrak{P}_i$ of $L$ lying above the prime $p_i$ such that 
$$\left(\frac{\mathfrak{P}_i}{L/\Q}\right)=\tau_i\; \text{for  each}\;  i.$$ 
Clearly $f(\mathfrak{P}_i|p_i)= \ell^{n_i}f$ for each $i$. Suppose $\mathfrak{p}_i=\mathfrak{P}_i\cap \mathcal{O}_K$, then $$f(\mathfrak{P}_i|p_i)=f(\mathfrak{P}_i|\mathfrak{p}_i) \times f(\mathfrak{p}_i|p_i) .$$ 
As $\ell$ is relatively prime to $n$ and $f(\mathfrak{p}_i|p_i)|n$, we conclude that 
$$f(\mathfrak{P}_i|\mathfrak{p}_i)=\ell^{n_i} \mbox{ and } f(\mathfrak{p}_i|p_i)=f.$$
If $$\left(\frac{\mathfrak{p}_i}{L/K}\right)=\tau'_i\; \text{for  each}\;  i.$$ Then we have $\tau_i^f=\tau'_i$ and order of $\tau'_i=\ell^{n_i}$ for each $i$. Observe that $\tau'_i$'s generate $G_{L/K}$. Thus, by the Artin isomorphism theorem the ideal classes of $\mathfrak{p}_i$'s generate the group $H$. Since $H$ contains $C\ell(K)[\ell]$, it shows that $f\in \mathbf{R}_{K/\Q}[\ell]$.\\

Let $f'$ be a positive divisor of $f$ and $f=af'$. Working with $\tau'=\tau^a$ in the above proof shows that $f'\in \mathbf{R}_{K/\Q}[\ell]$. 
\end{proof}

This completes the proof of Theorem \ref{A11} when $n=q$ is a prime. Similar arguments apply in the general case too. Next we prove Theorem \ref{A3}. 

\begin{proof} (Theorem \ref{A3})
Let $\sigma \in \Gal{\Q(\zeta_{\ell})/\Q}$ be an element of order $f$. As $f$ is relatively prime to $|Aut(C\ell(\Q(\zeta_{\ell})))|$, for any homomorphism $\psi : \Gal{Q(\zeta_{\ell})/\Q} \longrightarrow Aut(C\ell(\Q(\zeta_{\ell})))$ we must have $\psi(\sigma)=Id$. From Theorem \ref{Main}, we obtain $f \in \mathbf{R}_{\Q(\zeta_{\ell})/\Q}$. \\
Let $\theta_f$ be the annihilator as defined in (\ref{ann}). The Galois group $\Gal{\Q(\zeta_{\ell})/\Q}$ has unique subgroup of orders $f$ and $(\ell-1)/f$, say, $H_1$ and $H_2$ respectively.  By Lemma \ref{el}, we see that
\begin{equation}\label{51}
\theta_f =\sum_{\sigma \in H_2} \sigma.
\end{equation}
Let $F$ be the subfield of $\Q(\zeta_{\ell})$ of degree $f$. Then $\Gal{\Q(\zeta_{\ell})/F}=H_2$. Let $\mathfrak{P}$ be an unramified prime ideal of $\Q(\zeta_{\ell})$ and $\mathfrak{p}$ be the prime ideal of $F$ below $\mathfrak{P}$. From (\ref{51}) it follows that 
$$\theta_f(\mathfrak{P})=\mathfrak{p}^t \Z[\zeta_{\ell}],$$
where $t=f(\mathfrak{P}|\mathfrak{p})$ is a divisor of $(\ell -1)/f$. Since $\theta_f$ is an annihilator, the ideal $\mathfrak{p}^t \Z[\zeta_{\ell}]$ must be principal.  From the hypothesis, $t$ is relatively prime to $C\ell(\Q(\zeta_{\ell}))$. Hence $\mathfrak{p}\Z[\zeta_{\ell}]$ must be principal. Taking norm, we obtain
$$N_{\Q(\zeta_{\ell})/F}(\mathfrak{p}\Z[\zeta_{\ell}])=\mathfrak{p}^{(\ell-1)/f}$$ is principal.
Thus the exponent of the class group $H_F$ is at most $(\ell-1)/f$.
\end{proof}

Looking at the table of class numbers of maximal real cyclotomic fields in \cite{LW91}, it looks like that the class numbers $C\ell(\Q(\zeta_{\ell}))^{+}$ are a power of $2$ with good frequency. Our Theorem \ref{A3} can be used to conclude that, in many of such cases, the class group $$C\ell(\Q(\zeta_{\ell}))^{+} \cong \left(\frac{\Z}{2\Z}\right)^s,$$
for some integer $s$. Here, by many of the cases we mean when $(\ell-1)/2$ is odd, $C\ell(\Q(\zeta_{\ell}))^{-}$ is odd and $(\ell-1)/2$ is relatively prime to $\phi(C\ell(\Q(\zeta_{\ell}))^{-})$.\\


In the remaining part of this section we give two more consequences of our study. There is good overlap between proof of Theorem \ref{A3} and the following theorem, still we provide details.
\begin{thm}\label{A2}
Let $K/\Q$ be a cyclic extension of degree $6$ and $F$ be the intermediate field of degree $3$. Assume that $h_K$ is relatively prime to $6$, and $3$ does not divide $|Aut(C\ell(K))|$. Then the exponent of the class group $H_F$ is at most $2$. Moreover, there exists an integer $s \geq 0$ such that $H_F \cong (\Z/2\Z)^s$. 
\end{thm}

\begin{proof} (Theorem \ref{A2}) Using Theorem \ref{CMain}, we see that $3 \in \mathbf{R}_{K/\Q}$.  Consider the corresponding annihilator $\theta_3$ as defined in (\ref{ann}). 

The Galois group $G_K$ has unique subgroups of orders $3$ and $2$, say, $H_1$ and $H_2$ respectively. We observe that $\Gal{K/F} \cong H_2$.

From Lemma \ref{el}, we see that $$\theta_3=\sum_{\sigma \in H_2} \sigma.$$
Let $\mathfrak{P}$ be a prime ideal of $K$ and $\mathfrak{p}$ denote the prime of $F$ below $\mathfrak{P}$. If $f(\mathfrak{P}|\mathfrak{p})=t$ then $t|2$, and from the Dedekind theorem we obtain
$$\theta_3(\mathfrak{P})=\mathfrak{p}^t \mathbb{O}_K.$$
Using Theorem \ref{AT}, we see that $\mathfrak{p}^t \mathbb{O}_K$ is principal. As $t|2$ and $2$ does not divide $h_K$ we conclude that $\mathfrak{p} \mathbb{O}_K$ is principal. Now, taking norm gives
$$N_{K/F}(\mathfrak{p} \mathbb{O}_K)=\mathfrak{p}^2$$ is principal. Thus for any prime $\mathfrak{p}$ of $F$ the ideal $\mathfrak{p}^2$ is principal. Hence the exponent of $H_F$ is a divisor of $2$. This proves the Theorem.
\end{proof}

\begin{thm}\label{A4}
Let $K$ be a number field with class number bigger than one. If $H(K)$ denotes the Hilbert class field of $K$ then the extension $H(K)/\Q$ is never cyclic.
\end{thm}
\begin{proof}
If possible, suppose $H(K)/\Q$ is cyclic and the Galois group $G_{H(K)}$ is generated by $\sigma$. From Theorem \ref{CDT}, there is a prime ideal $\mathfrak{P}$ of $H(K)$ such that
\begin{equation}\label{e61}
\left( \frac{\mathfrak{P}}{H(K)/\Q} \right) =\sigma.
\end{equation}
Let $\mathfrak{p}=\mathfrak{P} \cap \mathbb{O}_K$ and $p=\mathfrak{P} \cap \Z$. Then 
\begin{equation}\label{e62}
f(\mathfrak{P}|p)=[H(K):\Q] \mbox{ and }f(\mathfrak{p}|p)=n.
\end{equation}
From this it follows that the ideal $\mathfrak{p}$ is a principal ideal. On the other hand, by Lemma \ref{FR} we have
\begin{equation}\label{e63}
\left( \frac{\mathfrak{p}}{H(K)/K} \right)=\left( \frac{\mathfrak{P}}{H(K)/K} \right)=\left( \frac{\mathfrak{P}}{H(K)/\Q} \right)^{f(\mathfrak{p}|p)}=\sigma^{n}.
\end{equation}
Since $\sigma^{n}$ generates the Galois group $G_{H(K)/K}$, from Lemma \ref{HCF} it follows that the class group $C\ell(K)$ is generated by prime ideal $\mathfrak{p}$. This implies that $C\ell(K)$ is trivial. This contradiction establishes the Theorem.
\end{proof}
We end this section with the remark that the consequences mentioned here are by no means exhaustive. Some more consequences will appear in second author's thesis.

\vspace{5mm}
\section{Explicit Examples}
We illustrate few examples explicitly, mainly among cyclotomic fields.  There are many examples, we illustrate few explicitly and record some in a table. For the value of class numbers we used the tables given in \cite{LW91}\\
1. Consider $K=\Q(\zeta_{163}^{+})$. Then $h_K=4$, and from Proposition \ref{CAAG} we obtain $|Aut(C\ell(K))|$ divides $6$. Using Theorem \ref{main} we see that $3, 9, 27 \in \mathbf{R}_{\Q(\zeta_{163}^{+})/\Q}$.

\vspace{1mm}
\noindent
2. Consider $K=\Q(\zeta_{191}^{+})$.  Then $h_K=11$. Using Theorem \ref{CMain} we see that $19 \in \mathbf{R}_{\Q(\zeta_{191}^{+})/\Q}$.

\vspace{1mm}
\noindent
3.  Consider $K=\Q(\zeta_{313}^{+})$. Then $h_K=7$. Using Theorem \ref{CMain} we see that $13 \in \mathbf{R}_{\Q(\zeta_{313}^{+})/\Q}$.

\vspace{1mm}
\noindent
4.  Consider $K=\Q(\zeta_{457}^{+})$. Then $h_K=5$. Using Theorem \ref{CMain} we see that $3, 19, 57 \in \mathbf{R}_{\Q(\zeta_{457}^{+})/\Q}$.

\vspace{1mm}
\noindent
5. Consider $K=\Q(\zeta_{457}^{+})$. Then $h_K=5$. Using Theorem \ref{CMain} we see that $3, 19, 57 \in \mathbf{R}_{\Q(\zeta_{457}^{+})/\Q}$.

\vspace{1mm}
\noindent
6. Consider $K=\Q(\zeta_{547}^{+})$. Then $h_K=4$, and from Proposition \ref{CAAG} we obtain $|Aut(C\ell(K))|$ divides $6$. Using Theorem \ref{CMain} we see that $7, 13 \in \mathbf{R}_{\Q(\zeta_{547}^{+})/\Q}$.

\vspace{1mm}
\noindent
7. Take $K=\Q(\zeta_{1399}^{+})$. Then $h_K=4$, and from Proposition \ref{CAAG} we obtain $|Aut(C\ell(K))|$ divides $6$. Using Theorem \ref{CMain} we get  $233 \in \mathbf{R}_{\Q(\zeta_{1399}^{+})/\Q}$.

\vspace{1mm}
\noindent
8. Take $K=\Q(\zeta_{1459}^{+})$. Then $h_K=247=13\times 19$. Using Theorem \ref{main} we obtain $3, 9 \in \mathbf{R}_{\Q(\zeta_{1459}^{+})/\Q}$.

\vspace{1mm}
\noindent
9. Suppose $K=\Q(\zeta_{1699}^{+})$. Then $h_K=4$, and from Proposition \ref{CAAG} we find that $|Aut(C\ell(K))|$ divides $6$. Using Theorem \ref{CMain} we get $283 \in \mathbf{R}_{\Q(\zeta_{1699}^{+})/\Q}$.

\vspace{1mm}
\noindent
10. Take $K=\Q(\zeta_{1879}^{+})$. Then $h_K=4$, and from Proposition \ref{CAAG} we obtain $|Aut(C\ell(K))|$ divides $6$. Using Theorem \ref{CMain} we see that $313 \in \mathbf{R}_{\Q(\zeta_{1879}^{+})/\Q}$.

\vspace{1mm}
\noindent
11. The class number of $\Q(\zeta_{96})$ is $9$. From Proposition \ref{Splitting-4}, we see that (\ref{44}) holds true. The group $Aut(C\ell(K))$ has either $6$ elements or $48$ elements. In either case, number of elements of order $2^r$ for some $r$ is at most $16$. On the other hand, all the elements of $\Gal{\Q(\zeta_{96})/\Q}$ have order $2^r$ for some $r$. Thus there is an element of order $2$ in $\Gal{\Q(\zeta_{96})/\Q}$ which is in the kernel of any homomorphism $\psi : \Gal{\Q(\zeta_{96})/\Q} \longrightarrow Aut(C\ell(K))$. Now proceeding as done in the proof of Theorem \ref{Main}, we can show that $2 \in \mathbf{R}_{\Q(\zeta_{96})/\Q}$.

\vspace{1mm}
\noindent
In the table below, we list some primes $\ell$, corresponding class numbers $C\ell(\Q(\zeta_{\ell}^{+}))$ of real cyclotomic fields,  and elements in $\mathbf{R}_{\Q(\zeta_{\ell}^{+})/\Q}$ which we could obtain from our study.  In most of the examples recorded here, it turns out that the class number is prime (even though it was not necessary for application of our results). There are many more examples possible, we do not try to list all the examples. This table is given just as an illustration of the outcomes of our study. Many more examples will appear in second author's thesis.\\

\begin{center}
	\begin{longtabu} to \linewidth {| l | l | l | l|}
\toprule
				$\ell$ &  $|C\ell(\Q(\zeta_{\ell}^{+}))$ & $\phi(|C\ell(\Q(\zeta_{\ell}^{+}))|)$ & $\mathbf{R}_{\Q(\zeta_{\ell}^{+})/\Q}$ \\
\midrule
				631   &    11  & 10 & 3, 7, 9,21, 63 \\
				761    &     3  & 2 & 5, 19, 95                 \\
				821   &    11  & 10 & 41                 \\
				829   &    47  & 46 &3, 9 \\
		857   &    5  & 4 &107 \\
			953   &    71  & 70 &17 \\ 
				977    &    5  & 4 & 61 \\ 
				1063    &    13  & 12 & 59 \\ 
				1069    &    7  & 6 &89 \\ 
				1093    &    5  & 4 &3, 7, 13, 21, 39, 91, 273 \\ 
				    1229 &    3  & 2 &307 \\ 
				1231    &    211  & 210 &41 \\ 
				1373    &    3  & 2 &7, 49, 343 \\ 
				1381    &    7  & 6 &5, 23, 115 \\ 
				1429  &    5 & 4 &3, 7, 17, 21, 51, 119,357 \\ 
				1567    &    7  & 6 & 29 \\ 
				1601    &    7 & 6 &5, 25 \\ 
				1697   &    17  & 16 & 53 \\ 
				1831    &    7  & 6 &5, 61, 305 \\ 
				1861   &    11 & 10 &3, 31, 93 \\ 
				1901   &    3  & 2 &5, 19, 25, 95, 475 \\ 
				1987    &    7  & 6 &331 \\ 
				2029   &    7  & 6 &13, 169 \\ 
				2113    &    37  & 36 &11 \\ 
				2213   &    3  & 2 &7, 79, 553 \\ 
				2351 &    11  & 10 &47 \\ 
				2381   &    11  & 10 &7, 17, 119 \\ 
				
				\hline
			
	\end{longtabu}
\end{center}


\vspace{5mm}
\section{Concluding remarks}
Computations of class numbers $C\ell(\Q(\zeta_{\ell}))^+$ of real cyclotomic field $\Q(\zeta_{\ell}^+)$ is an arduous task \cite{LF82, LW91}.  In a remarkable work \cite{RS03}, Schoof proposed a subgroup of $C\ell(\Q(\zeta_{\ell}))^+$ whose cardinality $\tilde{C\ell(\Q(\zeta_{\ell}))}$ can easily be computed. Conjecturally this subgroup equals to the class group $C\ell(\Q(\zeta_{\ell}^+))$. We wonder if our study can be used to lend credence to this conjecture. Using this conjectural equality we can obtain elements $f$ in $\mathbf{R}_{\Q(\zeta_{\ell}^+)/\Q}$. Now these $f's$ can be used to get information on class groups of intermediate fields. Many times, the class groups of these intermediate fields are relatively easier to study. A match between the actual computation and prediction from our studies should work as a support to the above mentioned conjecture of Schoof.\\

We believe that the study of primes of higher degree has good potential. It is desirable to obtain some analytic method to ensure some elements $f>1$ in $\mathbf{R}_{L/F}$ along the lines of $1\in \mathbf{R}_{L/F}$.

\end{document}